\documentclass[smallextended]{svjour3}       % onecolumn (second format)

\smartqed  % flush right qed marks, e.g. at end of proof
\usepackage{graphicx}
\usepackage{lineno, hyperref}
\usepackage[utf8]{inputenc}
\usepackage[T1]{fontenc}
\usepackage[draft=true]{minted}
\setminted[cpp]{fontsize=\footnotesize, frame=single}
\usepackage{amsmath}
\usepackage{amsfonts}
\usepackage{amssymb}
\usepackage{units}
\usepackage{soulutf8}
\usepackage{color}
\usepackage{fullpage}

\newcommand{\R}{\mathbb{R}}
\newcommand{\T}{\mathsf{T}}
\renewcommand{\L}{\mathcal{L}}
\renewcommand{\b}{\boldsymbol}

\newcommand{\x}{\b{x}}
\newcommand{\w}{\b{w}}
\newcommand{\eps}{\varepsilon}
\newcommand{\lap}{\nabla^2}

\newcommand{\einf}{e_\infty}
\newcommand{\Sph}[1]{B_{#1}}
\DeclareMathOperator{\Tr}{Tr}

\journalname{Journal of Scientific Computing}
\begin{document}

\title{Monomial augmentation guidelines for RBF-FD from accuracy vs.\ computational time perspective}


\author{Mitja Jan\v{c}i\v{c}         \and
	Jure Slak 			 \and
	Gregor Kosec
}


\institute{M. Jan\v{c}i\v{c} \at
	``Jo\v{z}ef Stefan'' International Postgraduate School, Jamova cesta 39, 1000 Ljubljana, Slovenia \\
	``Jo\v{z}ef Stefan'' Institute, E6, Parallel and Distributed Systems
	Laboratory, Jamova cesta 39, 1000 Ljubljana, Slovenia \\
	\email{Mitja.Jancic@ijs.si}           %  \\
	\and
	J. Slak \at
	``Jo\v{z}ef Stefan'' Institute, E6, Parallel and Distributed Systems
	Laboratory, Jamova cesta 39, 1000 Ljubljana, Slovenia \\
	Faculty of Mathematics and Physics, University of Ljubljana, Jadranska 19, 1000
	Ljubljana, Slovenia \\
	\email{Jure.Slak@ijs.si}           %  \\
	\and
	G. Kosec\at
	``Jo\v{z}ef Stefan'' Institute, E6, Parallel and Distributed Systems
	Laboratory, Jamova cesta 39, 1000 Ljubljana, Slovenia \\
	\email{Gregor.Kosec@ijs.si}
}

\date{Received: date / Accepted: date}

\maketitle

\begin{abstract}
	Local meshless methods using RBFs augmented with monomials have become increasingly
	popular, due to the fact that they can be used to solve PDEs on scattered node sets
	in a dimension-independent way, with the ability to easily control the order of the method,
	but at a greater cost to execution time.
	We analyze this ability on a Poisson problem with mixed boundary conditions
	in 1D, 2D and 3D, and reproduce theoretical convergence orders practically, also
	in a dimension-independent manner, as demonstrated with a solution of Poisson's
	equation in an irregular 4D domain.
	The results are further combined with
	theoretical complexity analyses and with conforming execution time measurements,
	into a study of accuracy vs.\ execution time trade-off for each dimension.
	Optimal regimes of order for given target accuracy ranges are extracted and presented,
	along with guidelines for generalization.
	\keywords{meshless methods \and RBF-FD \and Poisson's equation \and $n$-dimensional \and convergence rates \and optimal order}
\end{abstract}

\section{Introduction}
The Radial Basis Function-generated Finite Differences (RBF-FD), a local strong form mesh-free
method for solving partial differential equations (PDEs) that generalizes the traditional Finite
Difference Method (FDM), was first mentioned by Tolstykh~\cite{tolstykh2003using}.
Since then, the method has become increasingly popular~\cite{fornberg2015solving},
with recent uses in linear elasticity~\cite{slak2019refined}, contact
problems~\cite{slak2019adaptive}, geosciences~\cite{fornberg2015primer}, fluid
mechanics~\cite{kosec2018local}, dynamic thermal rating of power lines~\cite{maksic2019cooling},
advection-dominated problems~\cite{mavrivc2020equivalent,shankar2018hyperviscosity},
financial sector~\cite{milovanovic2018radial}, etc.

RBF-FD, similarly to other mesh-free methods, relies on approximation of differential
operators on scattered nodes, which is an important advantage over mesh-based methods, as
node generation is considered much easier than the mesh generation. In fact, mesh
generation is often the most cumbersome part of the solution procedure in traditional methods,
which, especially in 3D geometries, often requires significant assistance from the user. When meshless methods were first developed, many solutions used available mesh
generators for generating discretization nodes and discarding the connectivity information after the
mesh had been generated~\cite{liu2002mesh}. Such approach is computationally wasteful, does not
generalize to higher dimensions, and some authors even reported that it failed to generate distributions
of sufficient quality~\cite{shankar2018robust}.
Since then, various node positioning algorithms have been proposed. Popular algorithms
use iterative approaches~\cite{hardin2004discretizing,liu2010node},
advancing front methods~\cite{fornberg2015fast,lohner2004general} or sphere
packing methods~\cite{choi1999node}. In 2018, a pure meshless algorithm based
on Poisson disk sampling~\cite{bridson2007fast} was introduced. Later that year, the
first dimension-independent node generation algorithm that supported distributions with spatially
variable density appeared~\cite{slak2018generation}, where the authors also
demonstrated the stability of RBF-FD on scattered nodes, even for complex non-linear problems in
3D without any special treatment of stencil selection as proposed in~\cite{oanh2017adaptive}.
Instead, a cluster of nearest neighboring nodes proved to be a satisfactory stencil that can also
be efficiently implemented in dimension-independent code, using specialized data structures, such
as~$k$-d tree~\cite{yianilos1993data}.

A common drawback of often used RBFs, such as Gaussians or Hardy's multiquadrics,
is that they
include a shape parameter that crucially affects accuracy and stability of the
approximation~\cite{wendland2004scattered}. If the shape parameter is kept constant,
the method
converges, but stability issues arise when computing in the standard basis, due to
high condition
numbers of the collocation matrices. To fix the stability issue, more sophisticated
algorithms can
be used, such as RBF-CP, RBF-QA, RBF-GA and others~\cite{wrigth2017stable}, but
such methods
sometimes add significant additional costs. A simpler solution for the stability
issue is to scale the
shape parameter so that the product of the shape parameter and the nodal spacing is
constant. However,
this can lead to local approximations that are not convergent - this phenomenon has been called
lack of convergence due to stagnation errors~\cite{flyer2016role}. Stagnation can be fixed 
by adding monomial terms that ensure consistency up to a certain order. This
technique has been used together with Polyharmonic splines (PHS) as RBFs, which have an
additional advantage of not having a shape parameter~\cite{bayona2017role}. In addition, the order
of added monomials directly effects the order of the RBF-FD approximation, effectively enabling
control over the convergence rate of the RBF-FD~\cite{bayona2019insight}. Various successful
applications of RBF-FD with PHS have since been demonstrated, both in 2D and in
3D~\cite{shankar2018robust,slak2018generation,bayona2017role}. The dimensional independence has
already been noted by, e.g.,\ Ahmad et al.~\cite{Ahmad2017}, but the high order RBF-FD has not yet
been thoroughly analyzed with computational efficiency in mind, as the authors were more focused on
solving the time-dependent part of the PDE of interest.

Although the RBF-FD formulation is dimension-independent, in the sense that the same formulation
can be used in 1D, 2D, 3D and higher,
translating this elegant mathematical formulation and algorithms into actual efficient computer
code is far from trivial. In this paper, we present a dimension-independent PDE solution procedure
based on our in-house dimension-agnostic implementation~\cite{medusa} of RBF-FD.
By dimension-agnostic implementation we refer to the fact that exactly the same
code can be used to solve problems in one, two, three or more dimensional spaces, while values of
parameters are optimised for each dimension separately. The paper describes all solution procedure
elements in detail and presents a thorough analysis of accuracy
and execution time in one, two and three dimensions, on a Poisson problem on scattered nodes
with mixed boundary conditions. To fully illustrate the dimension independence, a
solution of a
4-dimensional problem on an irregular domain is presented. A C++ implementation
of all discussed solution elements is freely available for download~\cite{code}.
%All examples in
%this paper were computed using the in-house open source Medusa library~\cite{medusa}.

The rest of the paper is organized as follows: In section~\ref{sec:sol-proc},
the RBF-FD solution
procedure is presented, in section~\ref{sec:pois}, the model problem is investigated,
in section~\ref{sec:ex}, an additional example is shown, and in section~\ref{sec:con}, the
conclusions are presented.

\section{RBF-FD solution procedure}
\label{sec:sol-proc}
In this section, the main steps of the RBF-FD solution procedure are described.
First, the domain is populated with scattered nodes. Once the nodes are positioned,
in each discretization node the approximation of the partial differential
operator is performed,
resulting in stencil weights. Finally, in the PDE discretization phase, the PDE is transformed
into a system of linear equations, whose solution stands for a numerical solution
of the considered PDE.

\subsection{Positioning of nodes}
\label{sec:nodes}
In the node generation algorithm, candidate nodes are generated on a $d$-sphere in
a $d$-dimensional space. This effectively means that the node positioning algorithm
remains the same for every number of dimensions. However, some parameters,
e.g.\ the number of candidates, can be optimized for various numbers of dimensions.

The node positioning algorithm takes as an input a domain $\Omega \subset \R^d$
with a spacing function $h\colon \Omega \to (0, \infty)$ and optionally a list
of arbitrary starting ``seed
nodes'' $X \subset \Omega$, often distributed along the boundary. It returns a set of
nodes that are suitable for strong-form discretizations and distributed over $\Omega$
with mutual spacing around a point $p$ approximately $h(p)$.

\begin{figure}
	\centering
	\includegraphics[height=5cm]{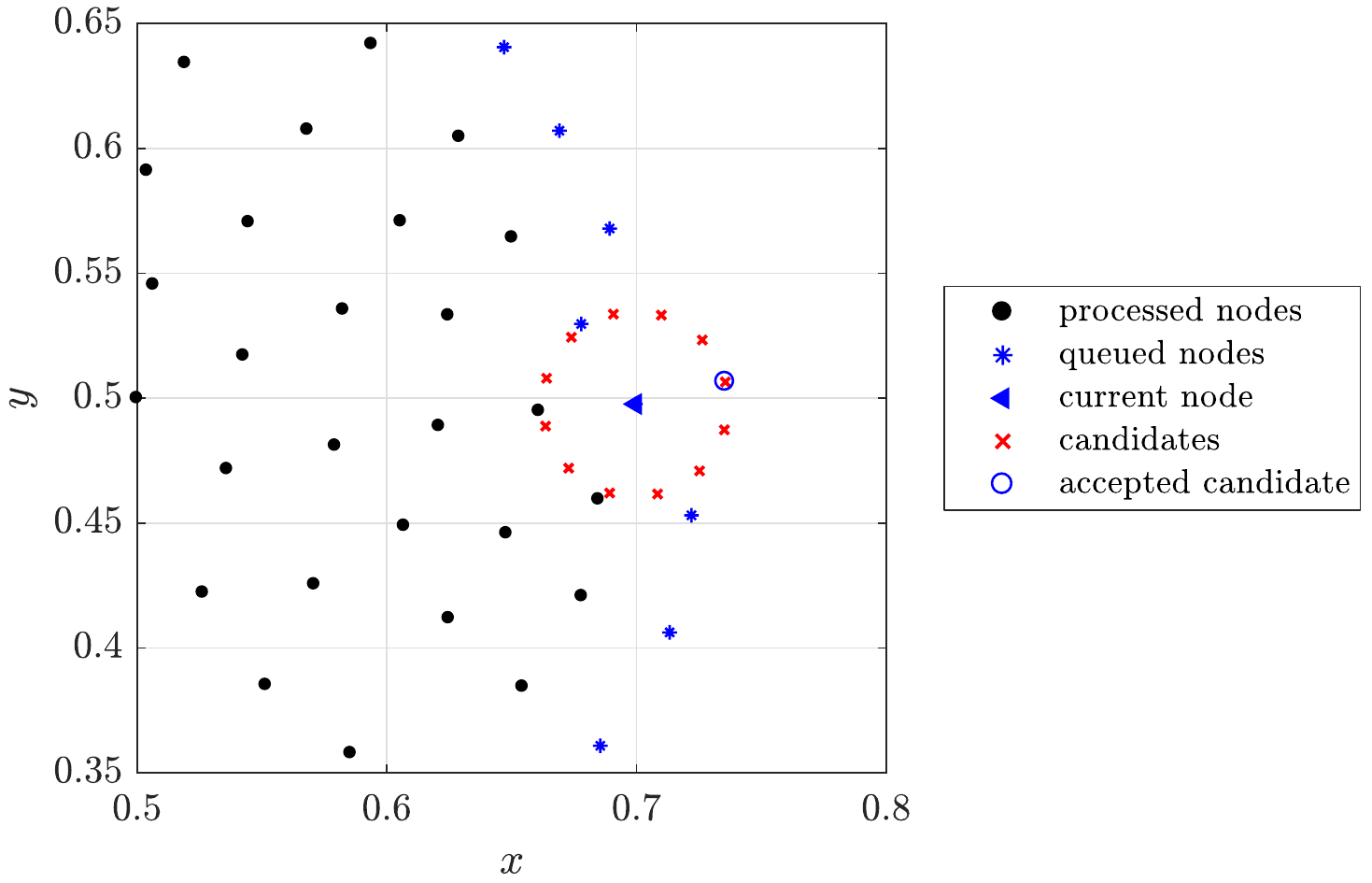}
	\includegraphics[height=5.5cm]{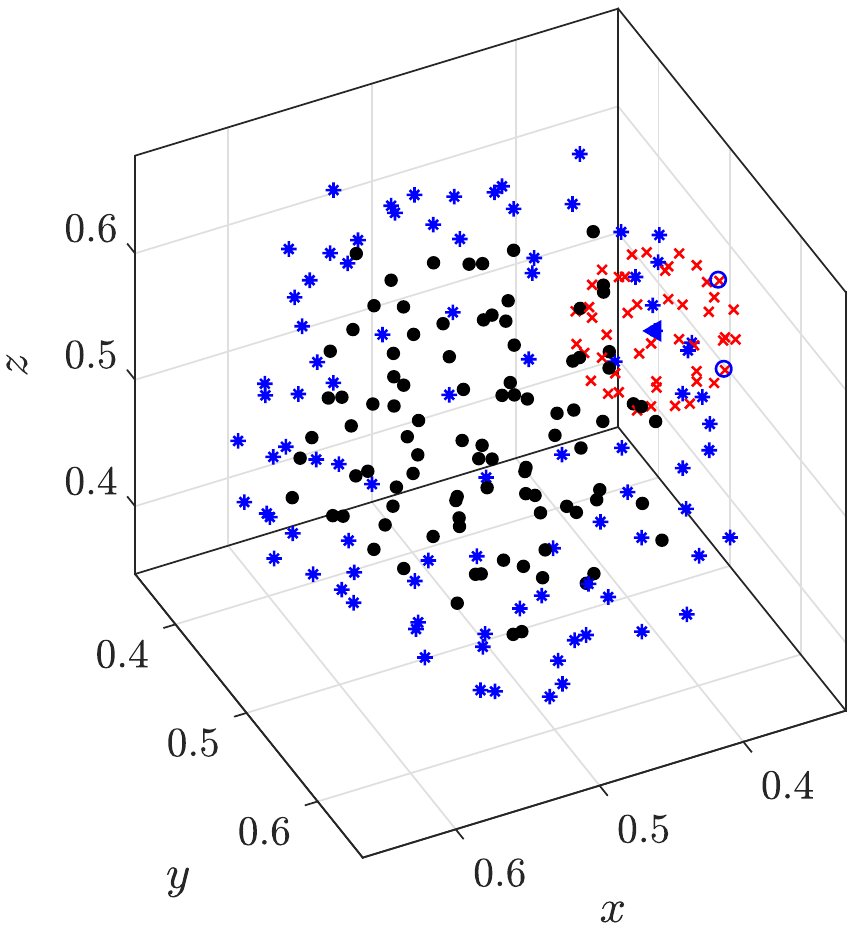}
	\caption{Node positioning algorithm during candidate generation phase.}
	\label{fig:fill}
\end{figure}

The algorithm used in this paper processes nodes in the input list in order.
For each node $p$, a number of expansion candidates distributed uniformly on a sphere centered at
$p$, of radius $h(p)$, are examined. If a candidate is inside the domain and sufficiently away from
the already processed nodes, it is accepted and added to the list $X$.
During the course of the algorithm, the list $X$ is implicitly partitioned into
already processed nodes, the current node, and future queued nodes.
Figure~\ref{fig:fill} shows this partition at a selected iteration in 2D and 3D, along
with the generated candidates from the current node, and flags the accepted ones.

Once all the elements of the list $X$ have been processed,
$X$ is returned as the resulting set of discretization nodes.
Further details and analyses of the algorithm are available in~\cite{slak2018generation}.
The stand-alone implementation of the algorithm is available online~\cite{standaloneFill} and
also included as a part of our in-house implementation of RBF-FD, the \emph{Medusa} library~\cite{medusa}.

\subsection{Approximation of partial differential operators}
\label{sec:rbffd}

Consider a partial differential operator $\L$
at a point $\x_c$. Approximation of $\L$ at a point $\x_c$
is sought using an ansatz
\begin{equation}
	\label{eq:approx}
	(\L u)(\x_c) \approx \sum_{i=1}^{n} w_i u(\x_i),
\end{equation}
where $\x_i$ are the neighboring nodes of $\x_c$ which
constitute its \emph{stencil}, $w_i$ are called \emph{stencil weights},
$n$ is the \emph{stencil size} and $u$ is an arbitrary function.

This form of approximation is desirable, since operator $\L$ at point $\x_c$
is approximated by a linear functional $\w_\L (\x_c)^\T$, assembled of weights $w_i$,
\begin{equation}
	\label{eq:approx-vec}
	\L|_{\x_c} \approx \w_\L (\x_c)^\T
\end{equation}
and the approximation is obtained using just a dot product with the function values
in neighboring nodes. The dependence of $\w_\L (\x_c)^\T$ on $\L$ and $\x_c$ is often omitted,
with $\w_\L (\x_c)^\T$ written simply as $\w$.

To determine the unknown weights $\w$, equality of~\eqref{eq:approx}
is enforced for a given set of basis functions. A natural choice are monomials, which
are also used in FDM, resulting in the Finite Point Method~\cite{onate1996finite}. However, using
monomial basis suffers from potential ill conditioning~\cite{mairhuber1956haar}. An alternative
approach is using an RBF basis.

In the RBF-FD discretization, the equality is satisfied for radial basis functions
$\phi_j$. These are RBFs, generated by a function $\phi\colon [0, \infty)
	\to \R$, centered at neighboring nodes of $\x_c$, given by
\begin{equation}
	\phi_j(\x) = \phi(\|\x-\x_j\|).
\end{equation}
Each $\phi_j$, for $j = 1, \ldots, n$, corresponds to one linear equation
\begin{equation}
	\sum_{i=1}^{n} w_i \phi_j (\x_i) = (\L \phi_j)(\x_c)
\end{equation}
for unknowns $w_i$. Assembling these $n$ equations into matrix form,
we obtain the following linear system:
\begin{equation} \label{eq:rbf-system}
	\begin{bmatrix}
		\phi(\|\x_1 - \x_1\|) & \cdots & \phi(\|\x_n - \x_1\|) \\
		\vdots                & \ddots & \vdots                \\
		\phi(\|\x_1 - \x_n\|) & \cdots & \phi(\|\x_n - \x_n\|)
	\end{bmatrix}
	\begin{bmatrix}
		w_1 \\ \vdots \\ w_n
	\end{bmatrix}
	=
	\begin{bmatrix}
		(\L\phi(\|\x-\x_1\|))|_{\x=\x_c} \\
		\vdots                           \\
		(\L\phi(\|\x-\x_n\|))|_{\x=\x_c} \\
	\end{bmatrix},
\end{equation}
where $\phi_j$ have been expanded for clarity.

The above system can be written more compactly as
\begin{equation} \label{eq:rbf-system-c}
	\b{\mathrm{A}} \w = \b \ell_\phi.
\end{equation}
The matrix $\b{\mathrm{A}}$ is symmetric, and for some basis functions $\phi$ even positive
definite~\cite{wendland2004scattered}.

Many commonly used RBFs, such as Hardy's multiquadrics or Gaussians,
depend on a shape parameter, which governs their shape and consequently affects
the accuracy and stability of the approximation. In this work,
we use polyharmonic splines (PHS), defined as
\begin{equation}
	\phi(r) =
	\begin{cases}
		r^k,       & k \text{ odd}  \\
		r^k\log r, & k \text{ even}
	\end{cases},
\end{equation}
to eliminate the need for a shape parameter tuning where $r$ denotes the Euclidean distance
between two nodes. Without monomial augmentation, local approximations using only PHS are not
convergent, nor do we have any guarantees of solvability. However, if the approximation given
by~\eqref{eq:rbf-system} is augmented with polynomials, we obtain convergence and conditional
positive definiteness, provided that the stencil nodes form a polynomially unisolvent
set~\cite{wendland2004scattered}. Augmentation is
performed as follows: Let $p_1, \ldots, p_s$ be polynomials forming the basis of the space of
$d$-dimensional multivariate polynomials up to and including total degree $m$, with $s =
	\binom{m+d}{d}$. In addition to the RBF part of the approximation, an exactness
	constraint
\begin{equation} \label{eq:p-constraint}
	\sum_{i=1}^{s} w_i p_j (\x_i) = (\L p_j)(\x_c)
\end{equation}
for monomials, is enforced. These additional constraints make the approximation overdetermined,
which is treated as a constrained optimization problem~\cite{flyer2016role}:
\begin{equation}
	\min_{\w} \left(\frac{1}{2} \w^\T \b{\mathrm{A}} \w - \w^\T \b \ell_{\phi}\right), \text{ subject to }
	\b{\mathrm{P}}^\T \b
	w = \ell_{p}.
\end{equation}
For practical computation, the optimal solution can be expressed as a solution of
a linear system
\begin{equation} \label{eq:rbf-system-aug}
	\begin{bmatrix}
		\b{\mathrm{A}}    & \b{\mathrm{P}} \\
		\b{\mathrm{P}}^\T & 0
	\end{bmatrix}\!\!\!
	\begin{bmatrix}
		\w \\ \b \lambda
	\end{bmatrix}
	=
	\begin{bmatrix}
		\b \ell_{\phi} \\ \b \ell_{p}
	\end{bmatrix}\!\!, \quad
	\b{\mathrm{P}} = \begin{bmatrix}
		p_1(\x_1) & \cdots & p_s(\x_1) \\
		\vdots    & \ddots & \vdots    \\
		p_1(\x_n) & \cdots & p_s(\x_n) \\
	\end{bmatrix}\!\!, \quad
	\b \ell_p = \begin{bmatrix}
		(\L p_1)|_{\x=\x_c} \\
		\vdots              \\
		(\L p_s)|_{\x=\x_c} \\
	\end{bmatrix},
\end{equation}
where $\b{\mathrm{P}}$ is a $n \times s$ matrix of polynomials evaluated at stencil nodes,
$\b \ell_p$ is the vector of values assembled by applying the considered operator $\L$ to
the polynomials at $\x_c$, and $\b \lambda$ are Lagrange multipliers.
Weights obtained by solving~\eqref{eq:rbf-system-aug}
are taken as approximations of $\L$ at $\x_c$, while values $\b \lambda$ are discarded.
The system~\eqref{eq:rbf-system-aug} is solvable if the stencil nodes form a polynomially
unisolvent set. This could potentially be problematic near the boundary, where it might happen that
all stencil nodes would be e.g.\ colinear or coplanar, but experience shows that this happens only
with stencil sizes which are too small to be a feasible approximation. With large enough stencil
sizes, stencils near the boundary always include at least some internal nodes. We did not use any
special techniques to ensure unisolvency, and did not run into any unisolvency-related issues.

The exactness of~\eqref{eq:p-constraint} ensures convergence behavior and control over the
convergence rate, since the local approximation has the same order as the polynomial basis
used~\cite{bayona2017role}, while the RBF part of the approximation~\eqref{eq:rbf-system}
takes care of potential ill-conditioning in purely polynomial approximation~\cite{flyer2016role}.
%to the Mairhuber--Curtis theorem~\cite{mairhuber1956haar}

\subsection{PDE discretization}
Consider the boundary value problem
\begin{align} \label{eq:bvp}
	\L u = f                  & \text{ in } \Omega,   \\
	u = g_d                   & \text{ on } \Gamma_d, \\
	\b n \cdot \nabla u = g_n & \text{ on } \Gamma_n,
\end{align}
with $\partial \Omega = \Gamma_d \, \cup \, \Gamma_n$, where the union is disjoint.
The domain $\Omega$ is discretized by placing $N$ scattered nodes $\x_i$
with quasi-uniform internodal spacing $h$, of which
$N_i$ are in the interior, $N_d$ on the Dirichlet and $N_n$ on the Neumann boundary.
Additionally, $N_{g}$ \emph{ghost} or \emph{fictitious} nodes are added outside the
domain on both Neumann and Dirichlet boundary, by translating the $N_d$ and the $N_n$ nodes on
$\partial \Omega$ for distance $h$ in the normal direction.

In the next step, stencils $\mathcal N(\b x_i)$ consisting of neighboring nodes
are selected for each node $\b x_i$. The most common approach is to compute stencils automatically,
by taking $n$ closest nodes for each node (including the node itself) as its stencil.

Next, partial differential operators appearing in the problem, such as $\L$ and $\partial_i$,
are approximated at nodes $\x_i$, using the procedure described in
section~\ref{sec:rbffd}. The computed stencils $\w_\L$ and $\w_{\partial_i}$
are stored for later use.

For each interior node $\x_i$, the equation $(\L u)(\x_i) = f(\x_i)$ is
approximated by a linear
equation
\begin{equation} \label{eq:disc-eq}
	\w_\L(\x_i)^\T \b u = \b f,
\end{equation}
where vectors $\b f$ and $\b u$ represent values of function $f$ and unknowns
$u$ in
stencil nodes of $\x_i$. For each Dirichlet boundary node $\x_i$, we have
the equation
\begin{equation} \label{eq:disc-dir}
	u_i = g_d(\x_i).
\end{equation}
For Neumann boundary nodes $\x_i$, the  linear equation
\begin{equation} \label{eq:disc-neu}
	\sum_{j=1}^d n_j \w_{\partial_j}(\x_i)^\T \b u = \b g_d
\end{equation}
approximates the boundary condition,
where similarly to before, vectors $\b g_d$ and $\b u$ represent values of
function
$g_d$ and unknowns $u$ in stencil nodes of $\x_i$.
Another set of $N_g$ equations is needed to determine the unknowns introduced
by ghost nodes.
Additionally to~\eqref{eq:disc-dir} and~\eqref{eq:disc-neu},
we also enforce~\eqref{eq:disc-eq} to hold for boundary nodes.

All $N_i + N_d + N_n + N_g$ equations are assembled into a sparse system with
$n(N_i + N_n + N_g) + N_d$ non-zero elements in general.
The solution $u_h$ of this system is a numerical approximation of $u$,
excluding the
values obtained in ghost nodes.

\subsection{Note on implementation}
We implemented the solution procedure described in this section in C++ using
object oriented approach and C++'s strong template system to achieve modularity
and consequent dimension independence.
The strongest advantage of the presented method is that all building blocks,
namely \emph{node positioning}, \emph{stencil selection}, \emph{differential operator
	approximation} and \emph{PDE discretization}, are independent and can therefore be
elegantly coded as abstract modules, not knowing about each other in the core
of their implementation. To ease the implementation of the solution procedure, additional
abstractions, such as \emph{operators}, \emph{basis functions}, \emph{domain shapes} and
\emph{approximations}, are introduced, acting as interfaces between the main blocks.
For example, to construct a RBF-FD approximation, one combines the RBF basis class with an augmented
\hbox{RBF-FD} class, computes stencil weights and supplies the computed weights into the
``operators'' class that enables the user to explicitly transform governing equations into the
C++ code, as demonstrated in the listing~\ref{lst:code}.

\begin{listing}
\begin{minted}{cpp}
// Define differential operator approximation.
Monomials<vec> mon(m);
Polyharmonic<double, k> ph;
RBFFD<decltype(ph), vec, ScaleToFarthest> appr(ph, mon);

// Compute stencil weights (shapes) with RBF-FD.
auto storage = domain.computeShapes<sh::lap|sh::d1>(appr);
Eigen::SparseMatrix<double, Eigen::RowMajor> M(N, N);
M.reserve(storage.supportSizes());
Eigen::VectorXd rhs(N); rhs.setZero();

// Prepare "operators" abstraction.
auto op = storage.implicitOperators(M, rhs);

// PDE discretization.
// Interior. 
for (int i : interior) {
    op.lap(i) = f_lap(domain.pos(i));
}
// Dirichlet boundary. 
for (int i : dir) {
    op.value(i) = f(domain.pos(i));
    op.lap(i, gh[i]) = f_lap(domain.pos(i));
}
// Neumann boundary.
for (int i : neu) {
    op.neumann(i, domain.normal(i)) = f_grad(domain.pos(i));
    op.lap(i, gh[i]) = f_lap(domain.pos(i));
}
\end{minted}
	\caption{A part of dimension-independent source code showing definition and
		sparse system assembly.}
	\label{lst:code}
\end{listing}

Vector and scalar fields are implemented as plain arrays, using a well developed linear algebra
library~\cite{eigenweb} that also implements or otherwise supports various direct and iterative
linear solvers. Please, refer to our open source Medusa library~\cite{medusa} for more
examples and features.

\section{Numerical example}
\label{sec:pois}

The behavior of the proposed solution procedure and its implementation are studied on a Poisson problem with
mixed boundary conditions. The aim is to analyze accuracy and convergence properties in one,
two and three dimensions.
%These convergence properties have been explored before, but never so
%systematically and only with Dirichlet boundary conditions~\cite{bayona2017role}.
Furthermore, theoretical computational complexity is discussed and
supported by experimental measurements of execution time, which
allows us to quantify the accuracy vs.\ execution time trade-off.

The problem is solved on an irregular domain $\Omega$, defined as
$\Omega = (\Sph 0 \cup \Sph 1 ) \setminus (\Sph 2 \cup \Sph 3)$, where
\begin{align}
	\Sph 0 & = \left \{\x \in \R^4, \ \left\|\b x - \b{\frac12}\right\| < \frac12 \right \},              \\
	\Sph 1 & = \left \{\x \in \R^4, \ \left\|\b x - \b{\frac 1 5}\right\|
	\le \frac 1 4 \right \},                                                                              \\
	\Sph 2 & = \left \{\x \in \R^4, \ \left\|\b x - \b{\frac12}\right\| \le \frac{1}{10} \right \} \text{
		and}                                                                                                  \\
	\Sph 3 & = \left \{\x \in \R^4, \ \left\|\b x - \b 1 \right\| \le \frac 1 2 \right \}
\end{align}
are balls in $\R^d$. For later use, the boundary $\partial\Omega$ is divided
into $\Gamma_d$ and $\Gamma_n$, the left and the right half
of the boundary, respectively
\begin{align}
	\Gamma_d & = \left \{ \x\in \partial \Omega, x_1 < \frac12 \right \} ,   \\
	\Gamma_n & = \left \{ \x\in \partial \Omega, x_1 \geq \frac12 \right \}.
\end{align}

\subsection{Governing equation}
\label{sec:def}
Numerical solution $u_h$ of Poisson's equation with both Dirichlet and Neumann
boundary conditions
is studied:
\begin{align}
	\label{eq:plap}
	\lap u(\x)   & = f_{lap}(\x) \qquad &                      & \text{in } \Omega,   \\
	\phantom{spaaaaaaaaaaaaaace}
	u(\x)        & = f(\x)              &                      & \text{on } \Gamma_d, \\
	\phantom{spaaaaaaaaaaaaaace}
	\nabla u(\x) & = \b f_{grad}(\x)
	\qquad       &                      & \text{on } \Gamma_n.
	\label{eq:pneu}
\end{align}
Here, the right hand side was chosen as
\begin{equation} \label{eq:closed_form}
	f(\x) = \frac{E(\x) }{g(\x)},
\end{equation}
where \begin{equation}
  E(\x) = \exp \Big( \sum _{i= 1}^d x_i^{a_i}\Big), \quad
  g(\x) = 1 + \x ^\T \b{\mathrm{H}} \x, \quad
  a_i  = 2 + i,  \quad
\end{equation}
$\b{\mathrm{H}}$ is a Hilbert matrix of size $d$, and $\hat{\b e}_i$ is the
$i$-th unit vector. The values of Laplacian and the gradient are computed from
$f$ as
\begin{align}
	f_{lap}(\x)     & = \frac{8 E(\x)}{g(\x)^3}
	(\b{\mathrm{H}}\x)^\T(\b{\mathrm{H}} \x) - \frac{2 E(\x)}{g(\x)^2}\bigg[2
	(\b{\mathrm{H}}\x)^\T (\sum_{i = 1}^d a_i x_i^{a_i - 1}\hat{\b e}_i) +
	\Tr(\b{\mathrm{H}})\bigg] \nonumber                                       \\
  & + \frac{E(\x)}{g(\x)} \bigg[\sum _{i = 1}^d a_i (a_i - 1) x_i^{a_i - 2} +
  (\sum _{i = 1}^d a_i x_i^{a_i - 1}\hat{\b e}_i)^\T (\sum_{i = 1}^d a_i
  x_i^{a_i - 1}\hat{\b e}_i)\bigg], \\
	\b f_{grad}(\x) & = \frac{E(\x)}{g(\x)} \bigg[\sum_{i = 1}^d a_i x_i^{a_i -
	1}\hat{\b e}_i - \frac{2}{g(\x)} (\b{\mathrm{H}} \x)^\T\bigg].
\end{align}
The closed-form solution $f$ of the above problem is a rational non-easily
separable function allowing us to validate the numerically obtained solution
$u_h$. The computed $u_h$ is only known at discretization points $\x_i$. The
errors between $u_h$ and $u$ are measured in three different norms:
\begin{align}
	e_1   & = \frac{\|u_h - u\|_1}{\|u\|_1}, \quad \|u\|_1 = \frac{1}{N} \sum_{i=1}^N |u_i|,
	\label{eq:e1}                                                                                     \\
	e_2   & = \frac{\|u_h - u\|_2}{\|u\|_2}, \quad \|u\|_2 = \sqrt{\frac{1}{N} \sum_{i=1}^N |u_i|^2},
	\label{eq:e2}                                                                                     \\
	\einf & = \frac{\|u_h - u\|_\infty}{\|u\|_\infty}, \quad \|u\|_\infty = \max_{i=1, \ldots, N}
	\label{eq:einf}
	|u_i|.
\end{align}

\begin{figure}[!h]
	\centering
	\includegraphics[width=\linewidth]{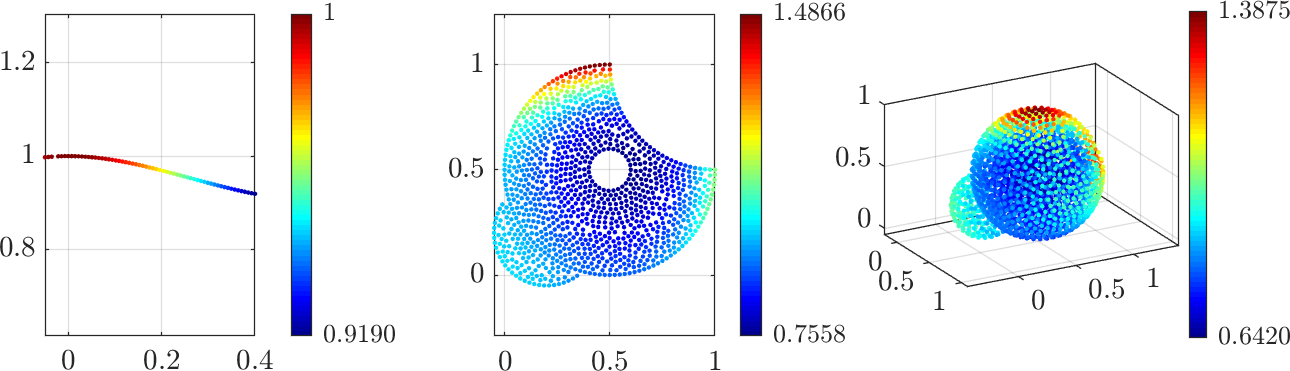}
	\caption{Computed numerical solution $u_h$ for $d=1, 2, 3$, from left to right. Chosen highest
		polynomial degree $m$ and node count $N$ are as follows: $N=64$ and $m = 4$ for $d=1$, $N=1286$
		and $m=2$ for $d = 2$ and $N = 3850$ and $m = 4$ for $d=3$.}
	\label{fig:solution}
\end{figure}

The problem~(\ref{eq:plap}--\ref{eq:pneu}) is studied in $d \in \{1, 2, 3\}$ dimensions.
Scattered computational nodes are generated using a dimension-agnostic node positioning
algorithm described in section~\ref{sec:nodes}. Ghost nodes were added to both Dirichlet
and Neumann boundaries, and are excluded from any post-processing. An example of
node distribution is shown in figure~\ref{fig:solution}.

Numerical results are computed using RBF-FD with PHS radial basis function
$\phi(r) = r^3$ and monomial augmentation, as described in
section~\ref{sec:sol-proc}. Radial function was kept same for all cases;
however, various orders of monomial augmentation were tested. For each dimension $d$,
solution to the problem is obtained using monomials up to and including degree
$m$, for $m \in \left \{-1, 0, 2, 4, 6, 8 \right \}$, where $m = -1$ represents a
pure RBF case with no monomials added. Only even orders of $m$ were used, because the
same order of convergence is observed with odd powers, but at a higher computational
cost~\cite{flyer2016role}.

Stencils for each node were selected by taking the closest $n$ nodes,
where $n$ was equal to two times the number of augmenting monomials, as
recommended by Bayona~\cite{bayona2017role}, or at least a FDM minimum of $2d+1$,~i.e.
\begin{equation}
	n = \max\left\{2\binom{m+d}{d},\ 2d+1\right\}.
\end{equation}
Specific values for $m$, $n$ and $d$ are presented in table~\ref{tab:ss}.

\begin{table}[h]
	\centering
	\renewcommand{\arraystretch}{1.2}
	\begin{tabular}{cccc} \hline
		\multicolumn{1}{c}{$m$} & \multicolumn{1}{c}{$d=1$} & \multicolumn{1}{c}{$d=2$} &
		\multicolumn{1}{c}{$d=3$}                                                             \\ \hline
		-1                      & 3                         & 5                         & 7   \\
		0                       & 3                         & 5                         & 7   \\
		2                       & 6                         & 12                        & 20  \\
		4                       & 10                        & 30                        & 70  \\
		6                       & 14                        & 56                        & 168 \\
		8                       & 18                        & 90                        & 330 \\ \hline
	\end{tabular}
	\caption{Support sizes in different dimensions for various augmentation orders.}
	\label{tab:ss}
\end{table}

BiCGSTAB with ILUT preconditioner was used to solve the sparse system.
Global tolerance was set to $10^{-15}$ with a maximum number of 500 iterations,
while the drop tolerance and fill-factor were dimension dependent:
$10^{-4}$ and $20$ for $d=1$, $10^{-4}$ and $30$ for $d=2$, and $10^{-5}$ and $50$ for $d = 3$,
respectively.

Figure~\ref{fig:solution} shows three examples of computed numerical solution $u_h$ for
each domain dimension $d$. The solutions are shown for various values of $m$ and
for small enough values of $N$ to also show nodal distributions.

\begin{figure}[!h]
	\centering
	\includegraphics[width=\linewidth]{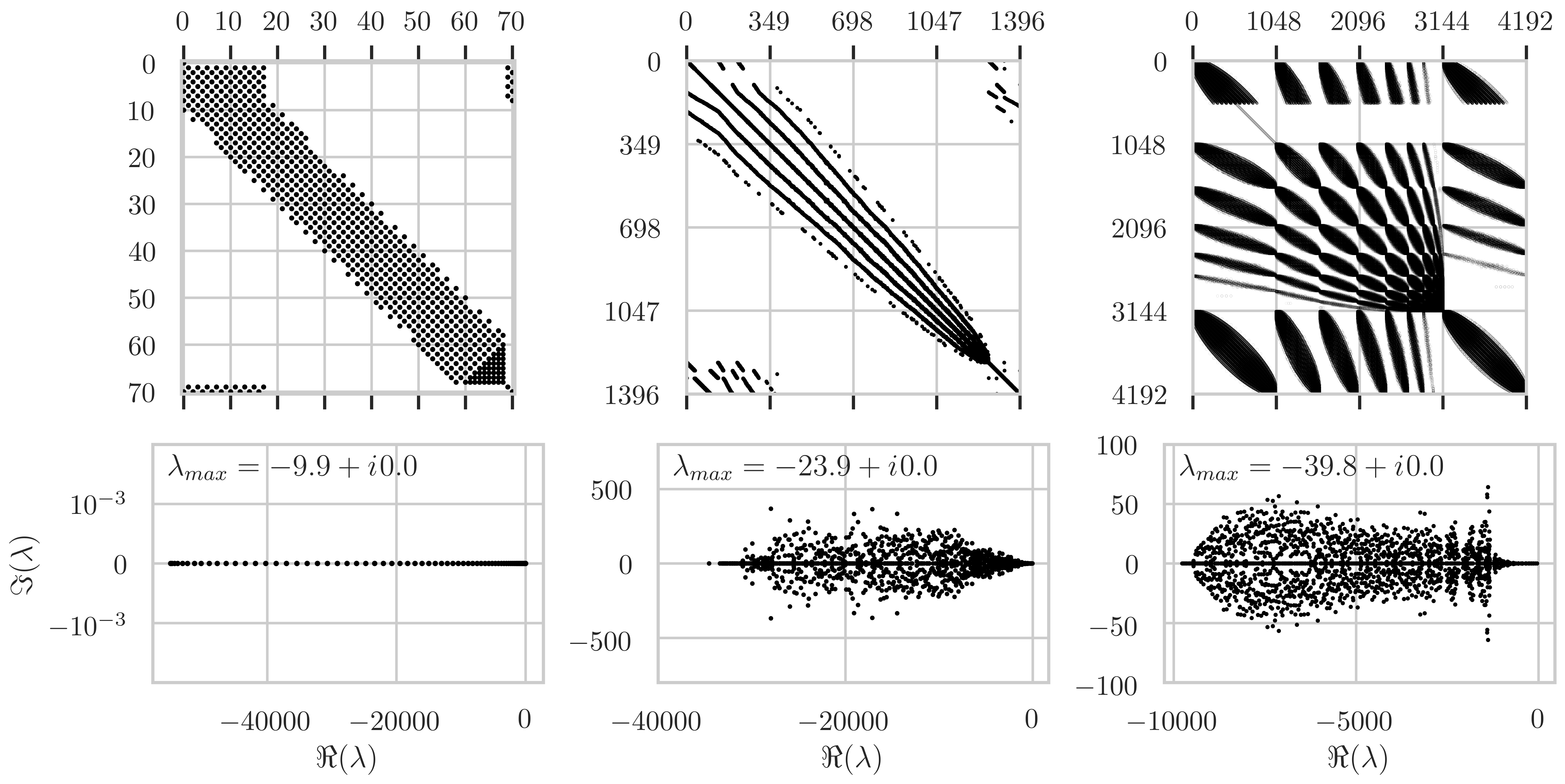}
	\caption{Plots of global sparse matrices (top row) and spectra of the Laplacian differentiation
		matrices (bottom row), corresponding to the solutions in figure~\ref{fig:solution}.}
	\label{fig:spectra}
\end{figure}

In the top row of figure~\ref{fig:spectra} global sparse matrices are shown. Additionally, spectra of the Laplacian differentiation matrices for cases shown in
figure~\ref{fig:solution} are shown in the bottom row of figure~\ref{fig:spectra},
to better assess the approximation quality. For all three cases, the eigenvalues have
negative real parts with relatively small spread around the imaginary axis.

\subsection{Convergence rate}

When using RBF-FD augmented with monomials, consistency is ensured up to order $m$,
which makes the expected convergence rate of at least $O(h^m)$. Here, $h$ denotes
the nodal spacing,
which is inversely proportional to $\sqrt[d]{N}$.

\begin{figure}[!h]
	\centering
	\includegraphics[width=\linewidth]{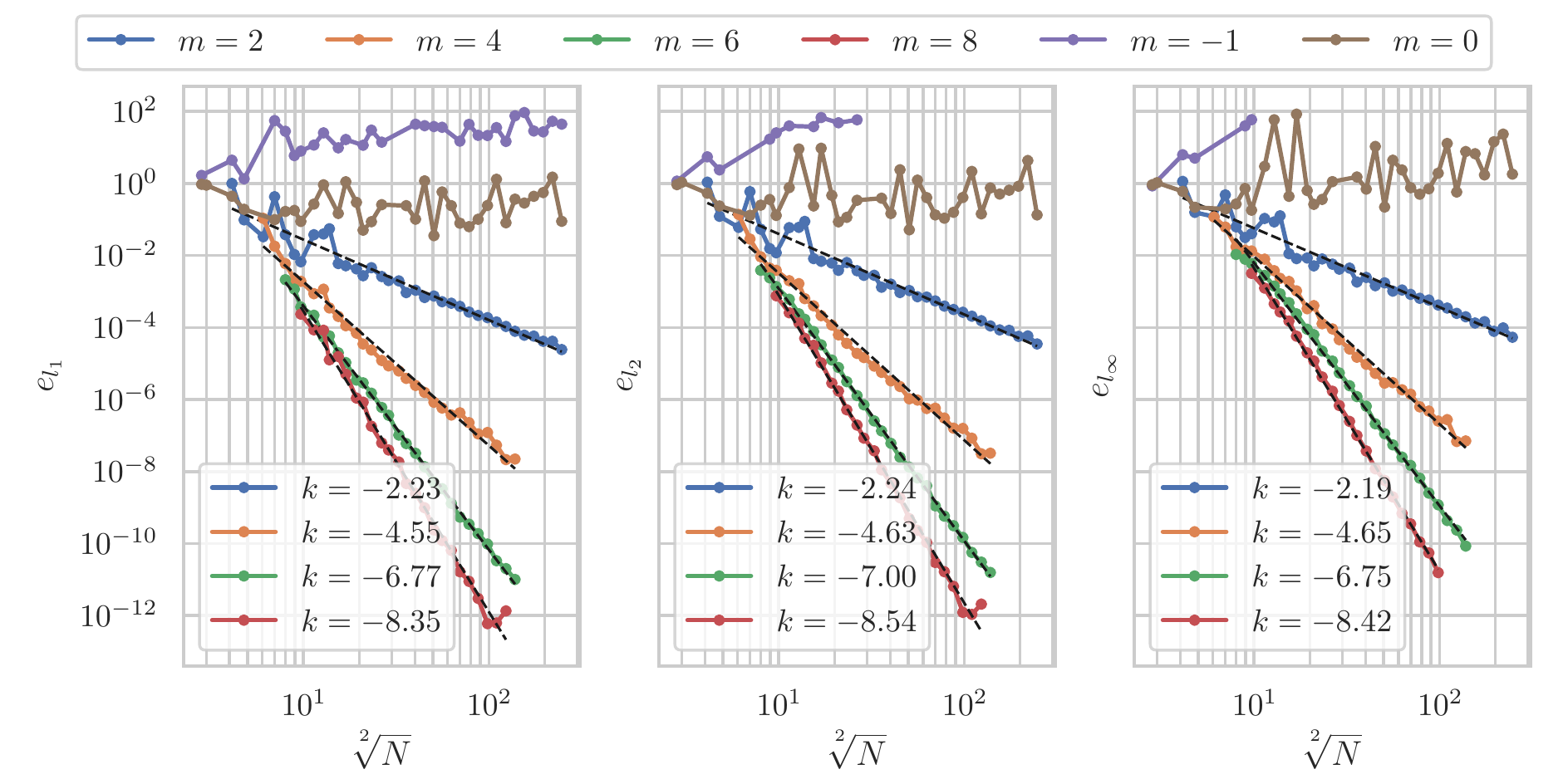}
	\caption{Errors between analytical solution $u$ and numerically obtained $u_h$,
		measured in three different norms. Computed are $e_1$, $e_2$ and $\einf$, from left to right,
		respectively, for the $d = 2$ dimensional case.}
	\label{fig:allerrors}
\end{figure}

Figure~\ref{fig:allerrors} shows $e_1$, $e_2$ and $e_\infty$ errors for various augmentation
orders in two dimensions.
The three errors have very similar values and similar convergence rates.
Convergence rates were estimated by computing the slope of a least-squares linear trend line
over the appropriate subset of the data.
Divergence is observed in the $m=0$ and the $m=-1$ case, which is consistent with properties of PHS
RBFs. These two cases are excluded from any further analyses in this paper.

In the rest of the discussion, only $\einf$ is used for convergence analysis,
since it measures the lowest convergence rates and
does not involve averaging, contrary to $e_1$ and $e_2$.

\begin{figure}[!h]
	\centering
	\includegraphics[width=\linewidth]{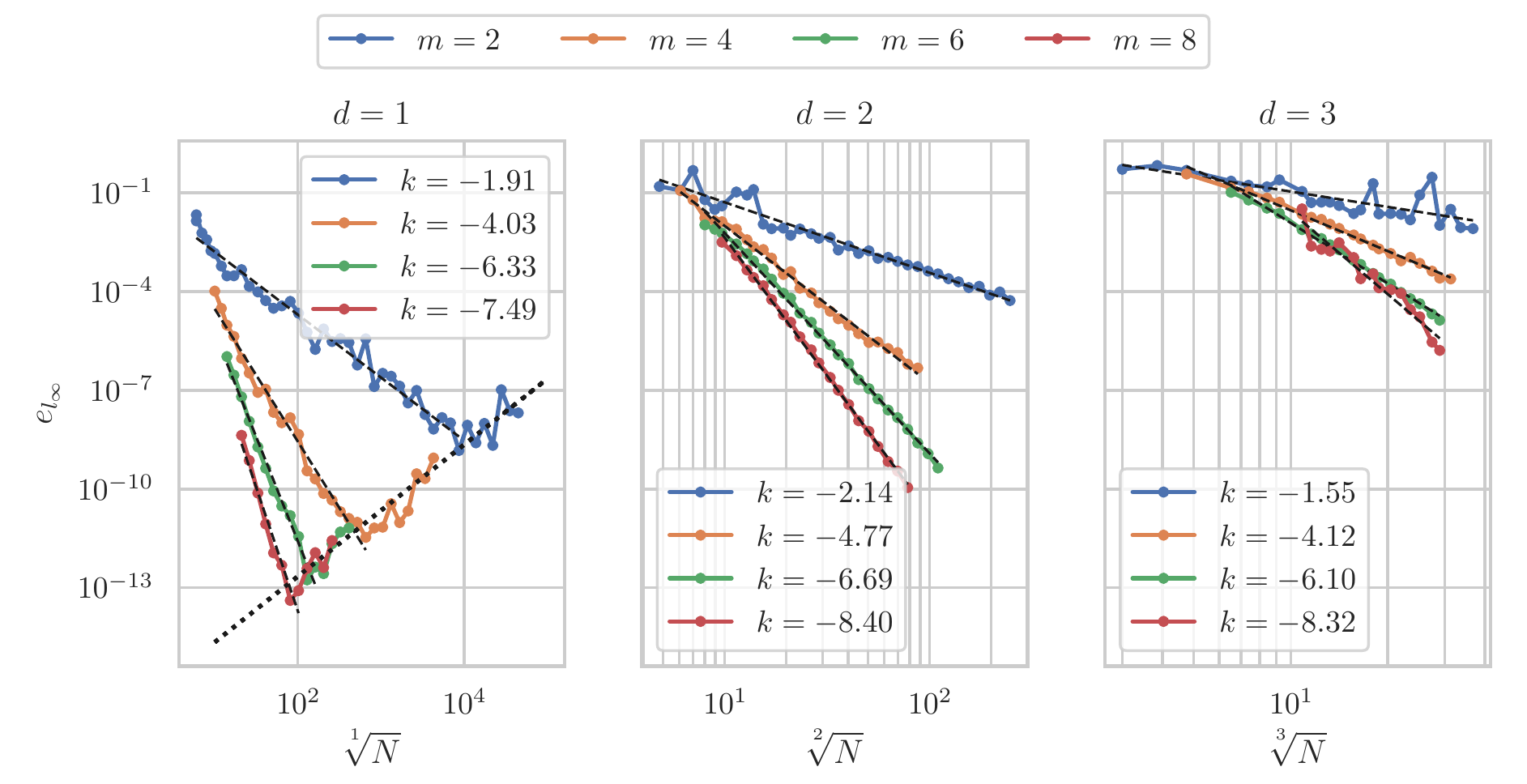}
	\caption{Convergence rate of $\einf$ for all domain dimensions $d = 1, 2, 3$, from
	left to right, respectively.}
	\label{fig:conv}
\end{figure}

Figure~\ref{fig:conv} shows the $\einf$ error for $d=1$, $d=2$, and $d=3$ dimensions.
The span of the horizontal axis was chosen in such a way that the total number of nodes in the
largest case was around $N = 10^5$ in all dimensions.
The observed convergence rates are independent of domain dimension and match the predicted order
$O(h^m)$. % they are also uncharacteristically smooth

All of the plots in the $d=1$ case eventually diverge, due to the errors in finite precision arithmetic, as previously noted for interpolation by Flyer et al.~\cite{flyer2016role}.
The dotted line in the $d=1$ case shows the $\eps/h^2$ line, where $\eps \approx 2.22 \cdot
	10^{-16}$.
The numerically obtained solution for the $d=3$ and $m=8$ case is unstable for smaller $N$.
For higher node counts $N$, the expected convergence behavior is obtained, as seen
from the fitted dashed line.

\subsection{Computational efficiency}
The importance of several different stages of $u_h$ computation is studied.
The computational procedure is divided into
\begin{itemize}
	\item \emph{node positioning}, where quasi-uniform placing of nodes in the domain
	      $\Omega$ and the domain boundary $\partial \Omega$, including positioning of $N_g$
	      ghost nodes, takes place. Node positioning time also includes finding the
	      stencils for each node in the domain,
	\item \emph{stencil weights computation}, where basis functions are defined and
	      shapes for the Laplace operator and first derivatives are stored,
	\item \emph{system assembly}, where computed weights are assembled in a sparse matrix and its
	      right-hand side is computed and
	\item \emph{system solution}, where the sparse system is solved.
\end{itemize}

\subsubsection{Computational complexity}
\label{sec:complexity}
The theoretical computational complexity is analyzed in this section. The total number of nodes
will be denoted as $N_t = N + N_g$; however, as $N_g$ nodes are distributed only along the
boundary, it holds that $N_g = O(N^{\frac{d-1}{d}})$ and thus $N_t = O(N)$.

The node positioning algorithm has complexity $O(N_t \log N_t)$~\cite{slak2018generation}. Finding
stencils of $n$ closest nodes takes $O(n N_t \log N_t)$ time, using a fast spatial search structure,
such as a $k$-d tree. The computation of stencils weights performs $N_t$ solutions of
linear systems of size $(n+s) \times (n+s)$, where $s = \binom{m+d}{d}$ is the number of
monomials used for augmentation. Since $n$ was chosen to be at least $2s$, it holds that $s = O(n)$.
Using LU decomposition or any other standard solution procedure for dense linear systems takes
$O((n+s)^3) = O(n^3)$ time. The total cost of weight computation is therefore $O(n^3N_t)$.

With appropriate pre-allocation of storage for the sparse matrix, system assembly takes linear
time in number of stencil nodes for each node, and right hand-side computation taken $O(1)$ per
node. The total cost of system assembly is thus $O(nN_t)$.

The solution of the sparse system uses iterative BiCGSTAB with ILUT preconditioner, whose
speed of convergence depends on the matrix properties.

The time complexity of the complete procedure is
\[
	O(nN_t\log N_t + n^3 N_t) + T,
\]
where $T$ is the complexity of the sparse solver.

\subsubsection{Execution time}
In this section, we measure execution time spent on different parts of the solution procedure.
All computations were performed on a single core of a computer with
\texttt{Intel(R) Xeon(R) CPU E5-2620 v3 @ 2.40GHz} processor and 64 GB of DDR4 memory.
Code was compiled using \texttt{g++ (GCC) 8.1.0} for Linux with \texttt{-O3 -DNDEBUG} flags.

\begin{figure}[!h]
	\centering
	\includegraphics[width=\linewidth]{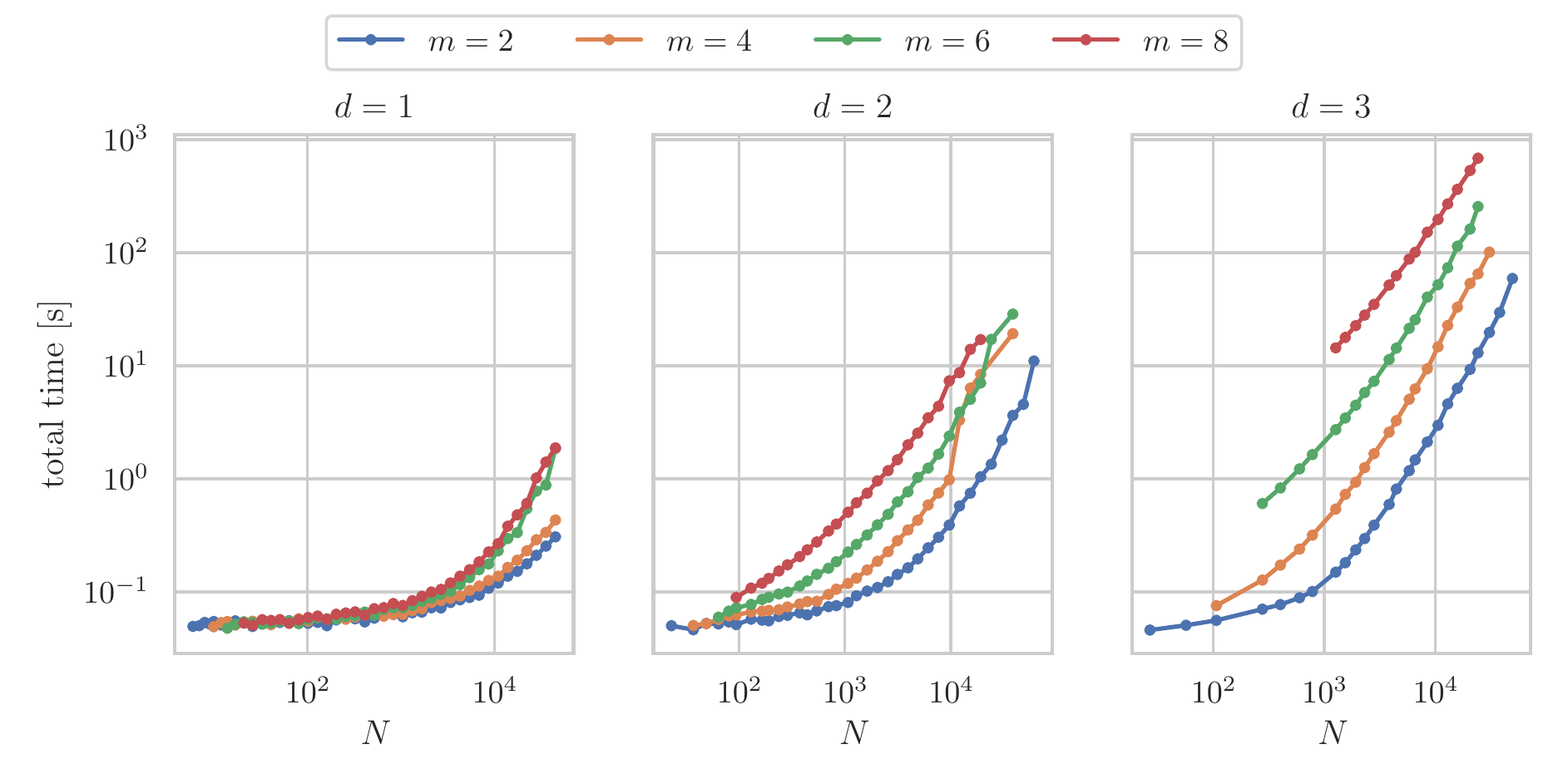}
	\caption{Median of 10 total execution times of $u_h$ computation for various setups.}
	\label{fig:totalTimes}
\end{figure}

Total execution times are shown in figure~\ref{fig:totalTimes} and
correspond to accuracy results in figure~\ref{fig:conv}.
The computational time grows with $N$ and with $m$, as expected from theoretical predictions
in section~\ref{sec:complexity}.

\begin{figure}[!h]
	\centering
	\includegraphics[width=\linewidth]{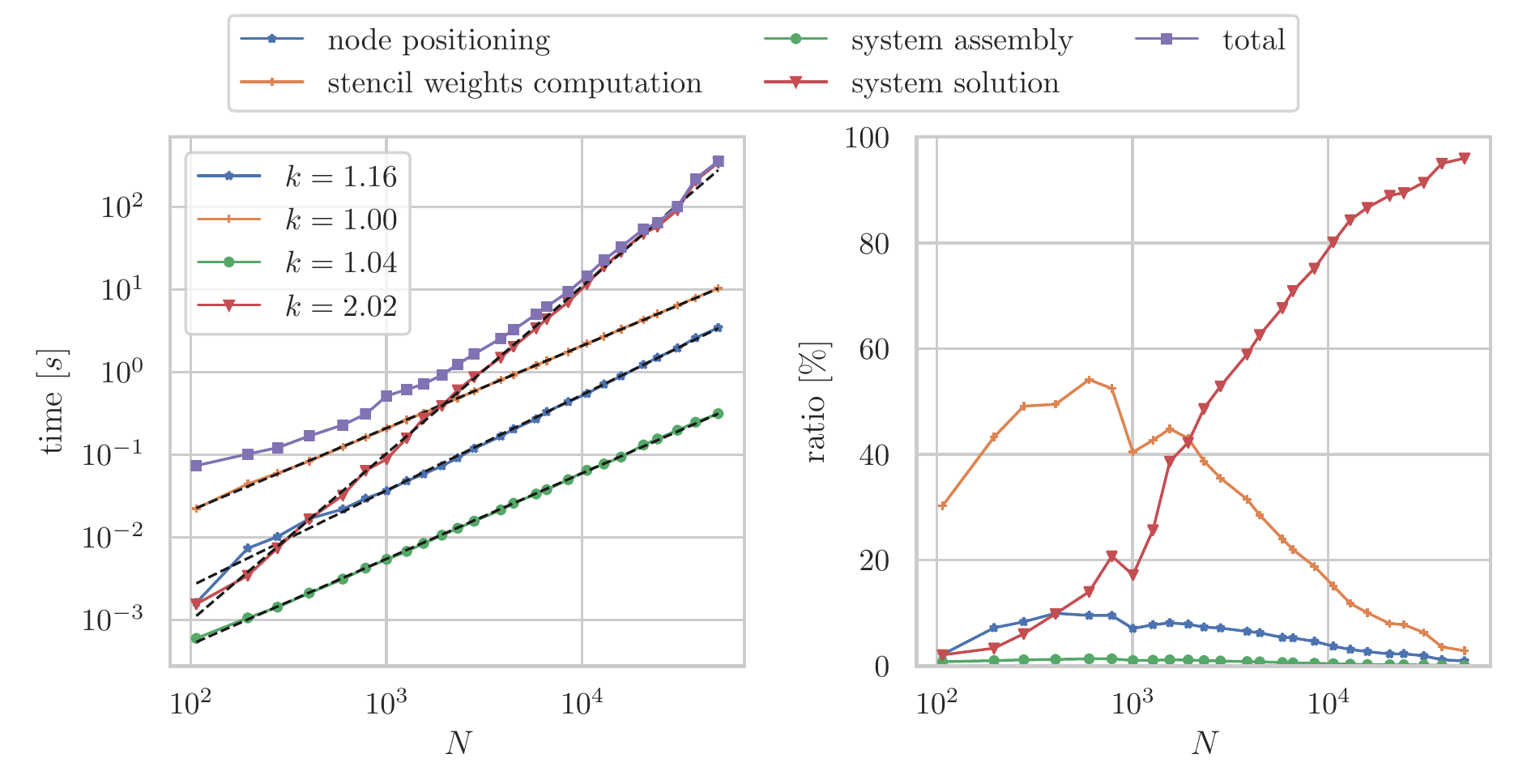}
	\caption{Absolute and relative times of different parts of the solution procedure
		for $d=3$ and $m = 4$.}
	\label{fig:times}
\end{figure}

Absolute times of different computation stages and their proportions to the total time are shown
in figure~\ref{fig:times}, on the left and the right side, respectively.
The observed growth rates match the theoretical complexities predicted for
node positioning, weight computation and system assembly.

Relative execution times provide additional insight into the execution of the solution
procedure and into optimization and parallelization opportunities.
The majority of the computational time is usually spent on either computing the stencil weights
(for smaller $N$) or on system solution (for large $N$).
Similar behavior was observed for other $m$ and in other dimensions, with
different percentage of total time spent on node positioning, weight computation and system
solution~\cite{kosec2019pareng}.

\subsection{Accuracy vs.\ execution time}
\begin{figure}[!h]
	\centering
	\includegraphics[width=\linewidth]{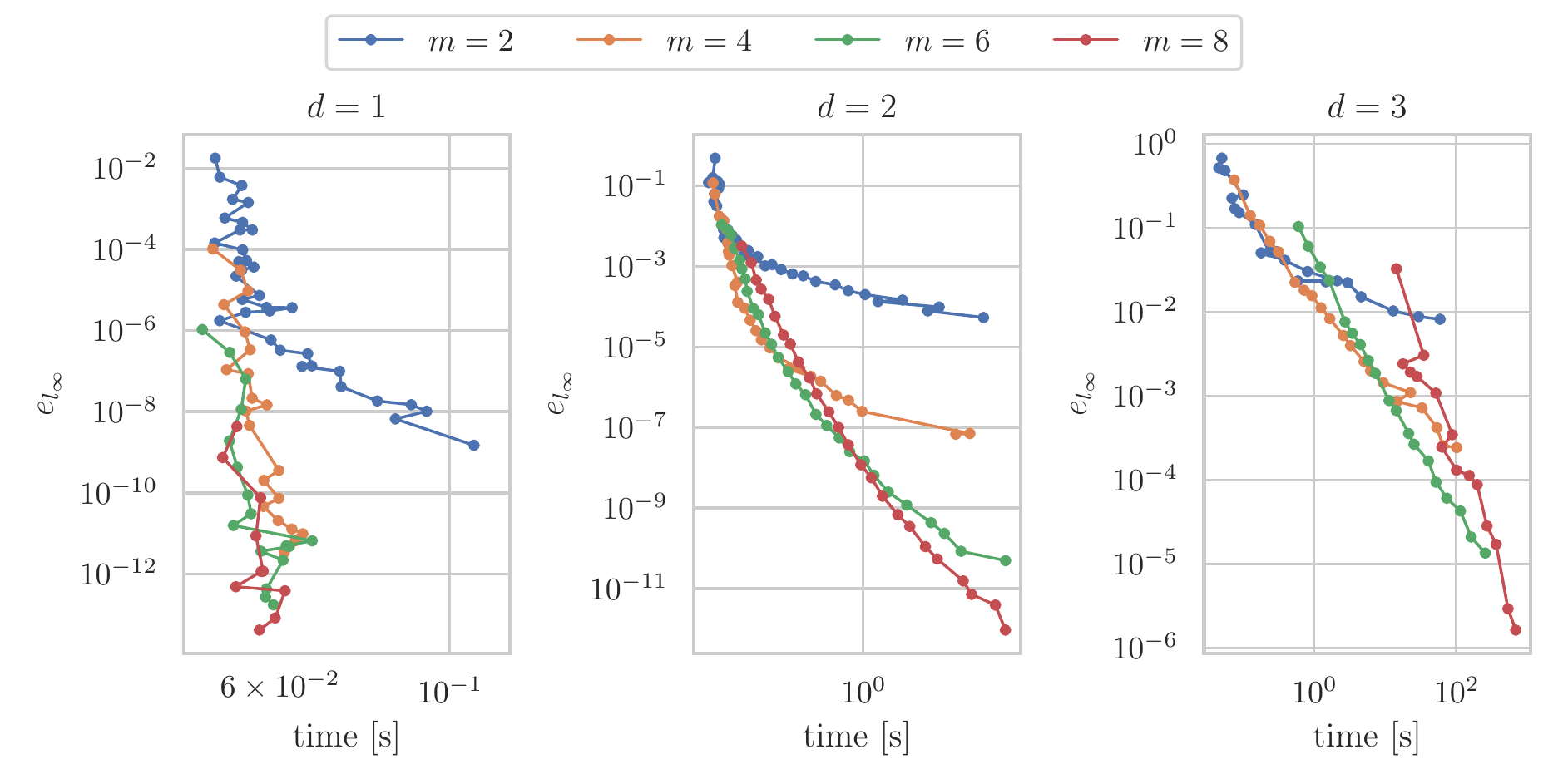}
	\caption{Accuracy vs.\ execution time trade-off for different orders of monomial augmentation.}
	\label{fig:efficiency}
\end{figure}

In the previous sections, we have shown that using higher orders, both accuracy
and execution time increase. In this section, we
analyze the accuracy vs.\ execution time trade-off.
Figure~\ref{fig:efficiency} shows $\einf$ error plotted with respect to the total computational time
needed to achieve it.

Significant differences can be observed between different orders of monomial augmentation.
For prototyping or any other sort of quick scanning of how or if the computed solution $u_h$
converges, using polynomials of a lower degree is undeniably very beneficial -- the computation of
$u_h$ takes little time, but at a cost of limited accuracy. When higher accuracy is required,
using polynomials of a higher degree can lead to a several orders faster computation time.
In some cases, using higher orders might even be a necessity,
e.g.\ for $d=2$, where accuracy of $\einf \approx 10^{-10}$ is reached the fastest by
$m=8$, while solution for $m=2$ would require $N$ out of reasonable computing capabilities.
The findings are summarized in table~\ref{tab:tradeoff}.

\begin{table}[h!]
	\centering
	\renewcommand{\arraystretch}{1.2}
	\begin{tabular}{cccccc} \hline
		\multicolumn{2}{c}{$d=1$}  &
		\multicolumn{2}{c}{$d=2$}  &
		\multicolumn{2}{c}{$d=3$}                                                                                        \\ \hline
		target accuracy $e_\infty$ & optimal $m$ & target accuracy $e_\infty$ & optimal $m$ & target
		accuracy $e_\infty$        & optimal $m$                                                                         \\
		\hline \rule{0pt}{12pt} % strut -- dont remove
		$10^{0}$ to $10^{-4}$      & 2           & $10^{0}$ to $10^{-2}$      & 2           & $10^{0}$ to $10^{-1}$  & 2 \\
		$10^{-4}$ to $10^{-6}$     & 4           & $10^{-2}$ to $10^{-5}$     & 4           & $10^{-1}$ to $10^{-3}$ & 4 \\
		$10^{-6}$ to $10^{-8}$     & 6           & $10^{-5}$ to $10^{-8}$     & 6           & $10^{-3}$ to $10^{-5}$ & 6 \\
		$10^{-8}$ to $10^{-13}$    & 8           & $10^{-8}$ to $10^{-12}$    & 8           & $10^{-5}$ to $10^{-7}$ & 8 \\
		\hline
	\end{tabular}
	\caption{Optimal setups for various desired target accuracy ranges in 1, 2 and 3 dimensions.}
	\label{tab:tradeoff}
\end{table}

Using the data in the table, we can extract a rough general recommendation.
As a rule of thumb, for the desired accuracy $e_\infty = 10^{-k}$ and dimension $d$, the
recommended order of augmentation is
\begin{equation} \label{eq:rot}
  m = \frac54 k + \frac45 d - 2,
\end{equation}
rounded to the nearest positive even integer. Even though the data points in the table are close to
being planar, the formula~\eqref{eq:rot} does not necessarily generalize well. A more general rule
is that the order of monomials should be increased with every two to three orders of increase in
accuracy, and that higher order augmentation should be more aggressively used in higher dimensions.

\section{Additional example}
\label{sec:ex}

In addition to already solved cases, we now demonstrate a solution of a 4-dimensional
Poisson problem~(\ref{eq:plap}--\ref{eq:pneu}). The irregular domain $\Omega$ is
now defined as
$\Omega = \Sph 0 \setminus (\Sph 1 \cup \Sph 2 \cup \Sph 3)$, where
\begin{align}
	\Sph 0 & = \left \{\x \in \R^4, \ \left\|\b x - \b{\frac12}\right\| < \frac12 \right \},          \\
	\Sph 1 & = \left \{\x \in \R^4, \ \left\|\b x - \left(\frac12, 1, \frac12, \frac12\right)\right\|
	\le \frac 1 4 \right \},                                                                          \\
	\Sph 2 & = \left \{\x \in \R^4, \ \left\|\b x - \b{0}\right\| \le \frac{13}{16} \right \} \text{
		and}                                                                                              \\
	\Sph 3 & = \left \{\x \in \R^4, \ \left\|\b x - \left(\frac12, \frac12, \frac 3 4,
	\frac12\right)\right\| \le \frac 1 8 \right \}
\end{align}
are balls in $\R^4$.

Dirichlet and Neumann boundary conditions are defined similarly to before, i.e.,\ $\Gamma_d$
is the
left half and $\Gamma_n$ is the right half of $\partial\Omega$. Additionally, the
boundary of the
smallest ball $\partial \Sph 3$ is added to the Dirichlet boundary:
\begin{align}
	\Gamma_d & = \left \{ \x\in \partial \Omega, x_1 < \frac12 \right \} \cup \partial \Sph 3,         \\
	\Gamma_n & = \left \{ \x\in \partial \Omega, x_1 \geq \frac12 \right \} \setminus \partial \Sph 3.
\end{align}

Scattered computational nodes were positioned using the same dimension-agnostic
node positioning
algorithm as before. A numerical solution $u_h$ was obtained using RBF-FD with PHS
$\phi (r) = r^3$ augmented with polynomials of degree $m=4$, according to
our rule of thumb~\eqref{eq:rot} for the desired accuracy $e_\infty = 10^{-2}$.

Approximately $N=85000$ nodes were positioned in $\Omega$ and closest $n=950$ nodes
were selected as
stencils for each node from the domain. Ghost nodes were, as in the previous case,
added to both Dirichlet and Neumann boundaries, and excluded from any post-processing.
The final system was solved using a direct sparse solver.

% BiCGSTAB with ILUT preconditioner was used to solve the sparse system, with max.\ iterations
% set to $2000$, global tolerance set to $10^{-15}$, drop tolerance set to $10^{-6}$ and fill
% factor set to $60$.

\begin{figure}[!h]
	\centering
	\includegraphics[width=0.8\linewidth]{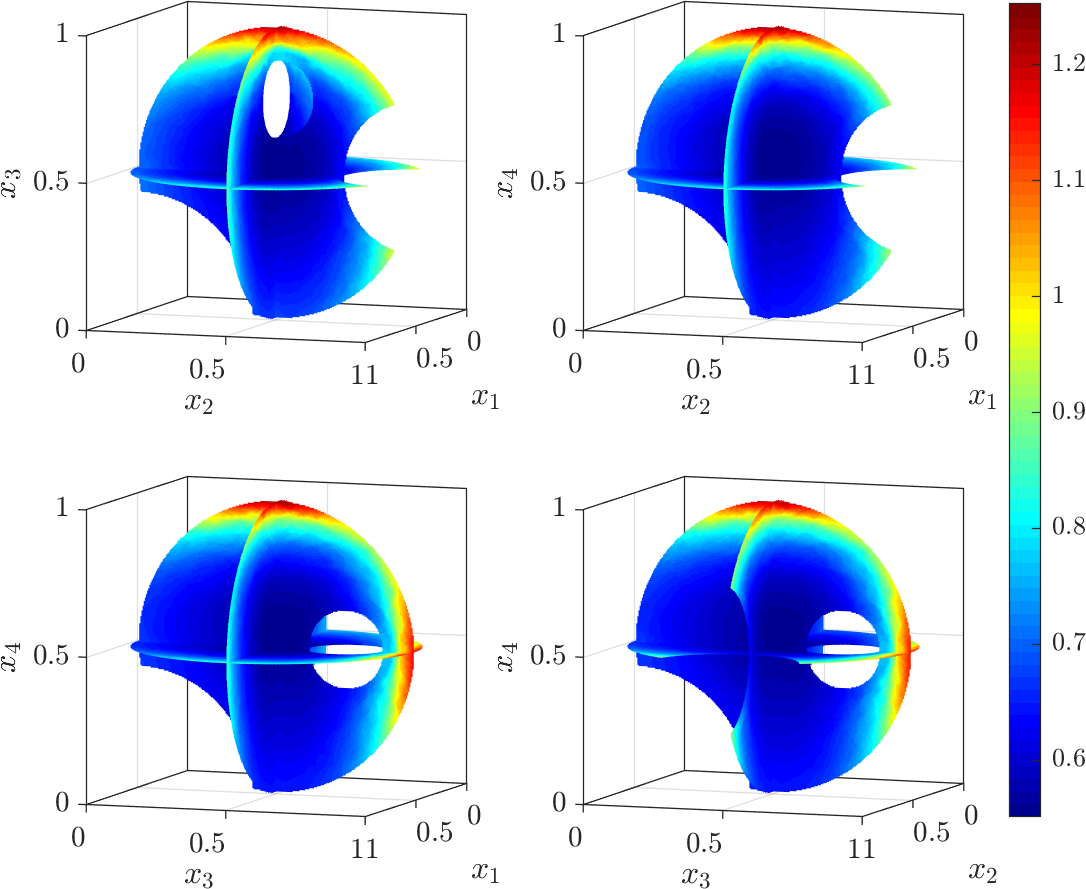}
	\caption{$3$-dimensional cross sections of a solution to a $4$-dimensional Poisson problem.}
	\label{fig:4Dsolution}
\end{figure}

Figure~\ref{fig:4Dsolution} shows the numerically obtained solutions. Four three-dimensional
slices are shown, defined by setting one coordinate to $x_i = 1/2$.
Modified Sheppard's interpolation algorithm~\cite{franke1980smooth} was used to interpolate the
solution to an intermediate grid, used for plotting the slices.

The solution is well-behaved even in 4 dimensions; however, a relatively large support size
is needed to obtain a desirable numerical stability. The errors equal to
$e_1 = 6.83\cdot 10^{-4}$, $e_2 = 2.11\cdot 10^{-3}$ and $e_\infty = 1.72 \cdot 10^{-2}$.
The total computational time spent was approximately \unit[15]{hours}.

\section{Conclusions}
\label{sec:con}
The message of this paper is twofold. First, we demonstrated that it is possible to design
an appropriately abstract implementation, which encompasses most of the meshless mathematical
elegance, allowing the user to construct a high order dimension-independent solution procedure.
To fully demonstrate the dimensional independence, we also presented a solution of a 4-dimensional
Poisson's problem on an irregular domain with both Neumann and Dirichlet boundary conditions.

Second, we used the devised implementation to analyze the increasing execution time that
comes tied with high order augmentation, to determine the conditions of
optimal computation efficiency for a desired target accuracy.

The analyses are performed on the solution of a Poisson problem with mixed boundary conditions in one,
two and three dimensions. To avoid shape parameter dependency, we used PHS augmented with monomials
as RBFs. Scattered nodes were positioned with a dedicated dimension-agnostic node generation
algorithm. The theoretical findings on how the highest order of the augmenting polynomial
directly controls the approximation rate of the RBF-FD independently of the domain dimension
are verified. A detailed breakdown of the computational complexity and the execution time
of different computational stages is also provided, to ensure that the implementation agrees with
the theoretical predictions. Finally, the high order vs.\ execution time trade-off is analyzed and the
findings are summarized in figure~\ref{fig:efficiency} and table~\ref{tab:tradeoff}.
While the analyses were done only for this particular problem, the results can be generalized in the
sense that for a high target accuracy, a high order method is a better choice, and vice versa.

Another interesting point are the increasing stencil sizes required for high
order methods, as shown in table~\ref{tab:ss}.
Especially in higher dimensions, this cost quickly becomes unmanageable.
%For example, in 3D
%for fourth order accuracy we already need 168 nodes in stencil.
Therefore, our future work will be
focused primarily on better understanding of the impact of the stencil size on the approximation quality.

%However, an interesting interplay of the stencil size, total number of nodes and the accuracy was discovered during the preparation of this paper
%domain dimension was discovered during preparation of this paper. Clarification of observed phenomena will be %the main focus of our future work.

\section*{Acknowledgements}
The authors would like to acknowledge the financial
support of the Slovenian Research Agency (ARRS) research core funding No.\ P2-0095
and the Young Researcher program PR-08346.

\bibliography{references}

\end{document}